\documentclass[12pt]{article}
\newlength{\mytopmargin}
\newlength{\myleftmargin}
\usepackage{amsmath,amsfonts,amssymb,amsthm}

\theoremstyle{plain}
\newtheorem{theorem}{Theorem}

\newtheorem{corollary}[theorem]{Corollary}
\newtheorem{prop}[theorem]{Proposition}

\theoremstyle{definition}

\theoremstyle{remark}

\numberwithin{equation}{section}

\begin{document}

\begin{center}
{\bfseries \LARGE On the gamma difference distribution  }\\[2\baselineskip]
{\Large Peter J. Forrester }\\[.5\baselineskip]

{\itshape 
School of Mathematics and Statistics, The University of Melbourne, Victoria 3010, Australia.\: \: Email: pjforr@unimelb.edu.au
}\\

\end{center}

\begin{abstract}
The gamma difference distribution is defined as the difference of two gamma distributions, with in general different shape and rate parameters. Starting with knowledge of the corresponding characteristic function, a second order linear differential equation characterisation of the probability density function is given. This is used to derive a Stein-type differential identity relating to the expectation with respect to the gamma difference distribution of a general twice differentiable function $g(x)$. Choosing $g(x) = x^k$ gives a second order recurrence  for the positive integer moments, which are also shown to permit evaluations in terms of ${}_2 F_1$ hypergeometric polynomials. A hypergeometric function evaluation is given for the 
absolute continuous moments. 
Specialising the gamma difference distribution gives the variance gamma distribution. Results of the type obtained herein have previously been obtained for this distribution, allowing for  comparisons to be made.
\end{abstract}

\section{Introduction}
\subsection{Functional form of the probability density function}
With $\alpha, \beta > 0$, and $\chi_A$ the indicator function for the condition $A$ ($\chi_A = 1$ for $A$ true and $\chi_A = 0$ otherwise) define the gamma random variable $X \in \Gamma[\alpha,\beta]$ by the probability density function
\begin{equation}\label{1.0}
p_X(x) = {\beta^\alpha \over \Gamma(\alpha)} x^{\alpha - 1} e^{-\beta x} \chi_{x > 0}.
\end{equation}
The corresponding characteristic function is
\begin{equation}\label{1.1}
\phi_X(t) = (1 - i t/\beta)^{-\alpha}.
\end{equation}
From the additive property of the characteristic function, it follows that with $X_1 \in \Gamma[\alpha_1,\beta_1]$ and $X_2 \in \Gamma[\alpha_2,\beta_2]$ we have
\begin{equation}\label{1.2}
\phi_{X_1 - X_2}(t) = (1 - i t/\beta_1)^{-\alpha_1}
(1 + i t/\beta_2)^{-\alpha_2}.
\end{equation}
The random variable $X_1 - X_2$ is said  to define the gamma difference distribution. The 2015 review of Klar \cite{Kl15} discusses a number of statistical properties and provides historical context. In the case that $X_1, X_2$ are identically distributed as chi-squared random variables, a comprehensive review is given in the recent work \cite{JO22}. Another recent work containing comprehensive review material is \cite{HGH22}, which furthermore goes into some depth in relation to the question of numerical evaluations. 

By Fourier inversion it follows from (\ref{1.2}) that the probability density function for the gamma difference distribution has the integral form
\begin{equation}\label{1.3}
p_{X_1 - X_2}(x) = \int_{-\infty}^\infty {e^{-i x t} \over (1 - i t/\beta_1)^{\alpha_1} (1 + i t/\beta_2)^{\alpha_2}} \, dt.
\end{equation}
An alternative integral form of the probability density function for $X_1 - X_2$ follows from the convolution structure
\begin{align}\label{1.5}
  p_{X_1 - X_2}(x) & = \int_{\max(0,-x)} p_{X_1}(x+x_2) p_{X_2}(x_2) \, dx_2 \nonumber \\
  & = {\beta_1^{\alpha_1} \beta_2^{\alpha_2} \over \Gamma(\alpha_1) \Gamma(\alpha_2)}
  \begin{cases}e^{\beta_2 x} 
  \int_x^\infty x_2^{\alpha_1 - 1} 
  (x_2 - x)^{\alpha_2 - 1}e^{- (\beta_1+\beta_2) x_2} \, dx_2, & x> 0 \\
 e^{-\beta_1 x} 
  \int_{-x}^\infty x_2^{\alpha_2 - 1} 
  (x_2 + x)^{\alpha_1 - 1}e^{- (\beta_1+\beta_2) x_2} \, dx_2, & x< 0.
  \end{cases}
\end{align}
It is noted in \cite[Eq.~(5)]{Kl15} that the integrals in (\ref{1.5}) can be found in \cite[Eq.~3.383(4)]{GR80}, implying an evaluation in terms of Whittaker's confluent hypergeometric function. It is furthermore the case that the Whittaker function can be replaced in favour of the Kummer hypergeometric function of the second kind (Tricomi function) $U(a,b;z)$ 
(also denoted $\Psi(a,;z)$) according to \cite[Eq.~(2.3)]{HGH22}
\begin{equation}\label{1.6}
  p_{X_1 - X_2}(x)  = {\beta_1^{\alpha_1} \beta_2^{\alpha_2} \over (\beta_1 + \beta_2)^{\alpha_1 + \alpha_2 - 1}}
  \begin{cases}
  {e^{-\beta_1 x} \over \Gamma(\alpha_1)}
 U(1-\alpha_1,2-\alpha_1 - \alpha_2; x (\beta_1 + \beta_2)),
   & x> 0 \\
{e^{\beta_2 x} \over \Gamma(\alpha_2)}
  U(1-\alpha_2,2-\alpha_1 - \alpha_1;- x (\beta_1 + \beta_2)),&  x < 0.
 \end{cases}
\end{equation}
This form can also be found in \cite{He17}. From the evaluation formula for $U(a,b;z)|_{z=0}$ with Re$(b) < 0$ \cite[Eq.~13.2.22]{DLMF} one notes from (\ref{1.6}) that
\begin{equation}\label{1.6a}
  p_{X_1 - X_2}(x)\Big |_{x=0}  = {\beta_1^{\alpha_1} \beta_2^{\alpha_2} \over (\beta_1 + \beta_2)^{\alpha_1 + \alpha_2 - 1}}
  {\Gamma(\alpha_1 + \alpha_2 - 1) \over \Gamma(\alpha_1) \Gamma(\alpha_2)}, \qquad \alpha_1 + \alpha_2 > 1;
  \end{equation}
  see \cite[displayed equation below (2.1)]{HGH22}.

  With $M(a,b;z)$ denoting the confluent hypergeometric function of the first kind (also denoted ${}_1 F_1(a,b;z)$), we know from \cite[Eq.~13.2.42]{DLMF} that
   \begin{equation}\label{pr.15} 
    U(a,b;z) = {\Gamma(1 - b) \over \Gamma(a-b+1)} M(a,b;z) + {\Gamma(b-1) \over \Gamma(a)} z^{1 - b} M(a-b+1,2-b;z).
 \end{equation}
Thus for $\alpha_1$ a positive integer
\begin{equation}\label{pr.16}
p_{X_1 - X_2}(x)  = {\beta_1^{\alpha_1} \beta_2^{\alpha_2} \over (\beta_1 + \beta_2)^{\alpha_1 + \alpha_2 - 1}}
  {\Gamma(\alpha_1 + \alpha_2 - 1) \over \Gamma(\alpha_1) \Gamma(\alpha_2)} e^{-\beta_1 x}
  M(1 - \alpha_1,2 - \alpha_1 - \alpha_2;x(\beta_1 + \beta_2)), \: \: x > 0,
  \end{equation}
  which has the analytic structure of an exponential times a polynomial. For more on this particular special case, see 
  \cite{Fo22} and \cite{KS20}.

Write $Y = X_1 - X_2$ for the special case of the gamma difference random variable $\alpha_1 =\alpha_2 = r/2$ and furthermore parameterised by setting $1/(\beta_1 \beta_2) = \sigma^2$, $(1/\beta_1 - 1/\beta_2) = 2 \theta$. Then from (\ref{1.3})
\begin{equation}\label{1.2a}
\phi_{Y}(t) = (1 - 2 i \theta t + \sigma^2 t^2)^{-r/2}.
\end{equation}
This characteristic function is better recognised not as specifying a particular gamma difference random variable, but rather as the characteristic function for the variance gamma distribution \cite[Eq.~(2.11)]{FGS23}. The latter is defined as the normal random variable N$[\theta a, \sigma^2 a]$ (mean $\theta a$, variance $\sigma^2 a$), in which $a$ itself is a random variable with distribution $\Gamma[r/2,1/2]$. The form given in the literature  of the corresponding probability density function does not involve the Tricomi function but rather the $K$-Bessel function (see \cite[Eq.~(1.1) with $\mu = 0$]{FGS23})
\begin{equation}\label{1.2b}
p_{Y}(x) = 
{e^{\theta x/\sigma^2} \over \sqrt{\pi \sigma^2}} {1 \over \Gamma(r/2)}
\Big ( {|x| \over 2 (\sigma^2 + \theta^2)^{1/2}} \Big )^{(r-1)/2}
K_{(r-1)/2} (|x| (\sigma^2 + \theta^2)^{1/2}/\sigma^2).
\end{equation}

A feature of $p_Y(x)$ known from \cite[Eq.~(2.26)]{FGS23} (see also \cite{Ga14}) is that it is the solution of a second order linear differential equation
 \begin{equation}\label{VG3} 
     {L}_x p_Y(x) = 0, \qquad   {L}_x := x  {d^2 \over d x^2} - \Big (  {2 \theta x \over \sigma^2} + (r-2) \Big ) {d \over dx} - \Big ( {x \over \sigma^2} - {(r-2) \theta \over \sigma^2} \Big ).
    \end{equation}  
    This can be used to show \cite{Ga14} that for a twice differentiable $g(x)$ and
    $Y=Y_{r,\theta,\sigma}$ a variance gamma random variable,
   \begin{equation}\label{rM}
   \mathbb E \Big ( \sigma^2 Y g''(Y) + (\sigma^2 r + 2 \theta Y) g'(Y) + (r \theta - Y) g(Y) \Big ) = 0,  
  \end{equation} 
  where it is assumed too that the growth of $g(Y)$ at $\pm \infty$ is such that all terms are well
  defined. Moreover, if (\ref{rM}) holds true, then $Y$ is a variance gamma random variable. Such an if and only if statement for a probability distribution satisfying a differential identity of the form (\ref{rM}) is known as a Stein characterisation; see \cite[\S 2.6]{FGS23}.

  We draw attention too to known formulas for the moments of $p_Y(x)$.
  In this regard, for $k$ a positive integer, set $\ell = \lceil k/2 \rceil + 1/2$ and $m = k\, {\rm mod} \, 2$. With $m_k(Y)$ denoting the $k$-th moment of the variance gamma distribution,
  $$
  m_k(Y) := \int_{-\infty}^\infty x^k p_Y(x) \, dx
  $$
  we have \cite[Eq.~(2.6)]{Ga22}, \cite[\S 2.7]{FGS23}
  \begin{equation}\label{2.6a}
  m_k(Y) = {2^{k+m} \theta^m \sigma^{r+2k} \over \sqrt{\pi} (\theta^2 + \sigma^2)^{(r+k+m)/2} \Gamma(r/2)}\Gamma\Big ({r-1 \over 2}+ \ell \Big ) \Gamma(\ell) \,
  {}_2 F_1 \Big ( \ell, {r-1 \over 2}+ \ell; {1 \over 2} + m; {\theta^2 \over \theta^2 + \sigma^2} \Big ).
  \end{equation}
  These moments satisfy the three term recurrence \cite[Eq.~(2.28)]{FGS23}
 \begin{equation}\label{2.6a+} 
 m_{k+1}(Y) = \theta (2k+r) m_k(Y) + \sigma^2 k (r+k-1) m_{k-1}(Y), \quad k \ge 0 ,
  \end{equation}
  with initial condition $m_0 = 1$. This implies that for $k$ even, $m_k(Y)$ is a polynomial in $r,\sigma,\theta$ each of degree $k$, and even in $\sigma, \theta$. For $k$ odd the recurrence implies $m_k(Y)$ is a polynomial in $r,\sigma,\theta$ of degree $k$ in $r,\theta$ and degree $k-1$ in $\sigma$, and is even in $\sigma$ and odd in $\theta$. Note that these structural properties are not immediate from the closed form (\ref{2.6a}).
  There is also a known closed form for the continuous absolute moment $\mathbb E(|Y|^k)$ with $k$ real and satisfying $k > {\rm max} \, \{-1,-r \}$. Thus from \cite[Eq.~(2.3)]{Ga22} and \cite[\S 2.7]{FGS23},
  \begin{equation}\label{2.6b}
  \mathbb E(|Y|^k) = {2^{k} \sigma^{r+2k} \over \sqrt{\pi} (\theta^2 + \sigma^2)^{(r+k)/2} \Gamma(r/2)}\Gamma\Big ({r+k \over 2} \Big ) \Gamma\Big ({k+1 \over 2} \Big ) \,
  {}_2 F_1 \Big ( {k+1 \over 2}, {r+k \over 2}; {1 \over 2} ; {\theta^2 \over \theta^2 + \sigma^2} \Big ).
  \end{equation}

  Our aims in this work are to give generalisations for the  gamma difference distribution of (\ref{VG3}), (\ref{rM}), (\ref{2.6a}), (\ref{2.6a+}) and (\ref{2.6b}).
  
\subsection{Summary of results}
Our first result begins with (\ref{1.3}) to deduce a differential equation for $p_{X_1 - X_2}(x)$. Associated with this is a Stein-type differential identity relating to the expectation of a general twice differentiable function $g(x)$ with respect to the gamma difference distribution.

\begin{prop}\label{P1}
For $p_{X_1 - X_2}(x)$ denoting the probability density function for the 
gamma difference distribution,
 \begin{multline}\label{2.7}
 x p_{X_1 - X_2}''(x) + \Big ( x (\beta_1 - \beta_2) + (2 - \alpha_1 - \alpha_2) \Big ) p_{X_1 - X_2}'(x) \\
 + \Big ( ( \beta_1 - \beta_2 + \alpha_1 \beta_2 - \alpha_2 \beta_1 ) - x \beta_1\beta_2 \Big ) p_{X_1 - X_2}(x) = 0.
 \end{multline}
 This is valid for general $\alpha_1, \alpha_2 > 0$ except at $x=0$ for $\alpha_1 + \alpha_2 < 1$ when $p_{X_1 - X_2}(x)$ diverges and so is not differentiable.
 Consequently, for $g(x)$ a twice differentiable function on $x \in \mathbb R$, and $X$ denoting a gamma difference random variable
 \begin{equation}\label{1.18}
 \mathbb E \Big ( X g''(X) + \Big ( X (\beta_1 - \beta_2) + (\alpha_1 + \alpha_2) \Big ) g'(X) +
(\alpha_1 \beta_2 - \alpha_2 \beta_1 - X \beta_1 \beta_2    ) g(X) \Big ) = 0,
\end{equation}
assuming the behaviour of $g(x)$ as $x \to \pm \infty$ is such that all implied expectations are well defined.
 \end{prop}

 A recurrence for the positive integer moments easily follows from (\ref{1.18}).

 \begin{corollary}\label{C1}
The $k$-th positive integer moment $m_k(X)$ of the gamma difference random variable $X$ satisfies the three term recurrence
\begin{equation}\label{1.19}
\beta_1 \beta_2 m_{k+1}(X) = \Big ( k (\beta_2 - \beta_1) + (\alpha_1 \beta_2 - \alpha_2 \beta_1) \Big )m_{k}(X) + k( k - 1 + \alpha_1 + \alpha_2) m_{k-1}(X), \quad k \ge 0,
\end{equation}
subject to the initial condition $m_0(X) =1 $.
 \end{corollary}

 We turn now to a closed form evaluation of the positive integer moments, and the absolute continuous moments.
 \begin{prop}\label{P3}
 The $k$-th positive integer moment $m_k(X)$ of the gamma difference random variable $X$ admits the ${}_2 F_1$ hypergeometric polynomial expressions
 \begin{align} \label{1.20}
m_k(X) & = {1 \over \beta_1^k}
\bigg ( \prod_{i=1}^k (\alpha_1 + i - 1) \bigg )
 \, {}_2 F_1 \Big ( -k, \alpha_2; - (\alpha_1 + k - 1); -{\beta_1 \over \beta_2} \Big ) \nonumber \\
 & = (-1)^k {1 \over \beta_2^k}
 \bigg ( \prod_{i=1}^k (\alpha_2 + i - 1) \bigg )
 \, {}_2 F_1 \Big ( -k, \alpha_1; - (\alpha_2 + k - 1); -{\beta_2 \over \beta_1} \Big ).
 \end{align}
 \end{prop}
 \begin{prop}\label{P4}
 Let $b > {\rm max} \,(0, 1 - \alpha_1 - \alpha_2)$.
 The $(b-1)$-th absolute continuous moment of the gamma difference random variable $X$ can be expressed in terms of the ${}_2 F_1$ hypergeometric function according to
 \begin{multline} \label{1.21} 
\mathbb E(|X|^{b-1}) = {\beta_1^{\alpha_1} \beta_2^{\alpha_2} \over (\beta_1 + \beta_2)^{\alpha_1 + \alpha_2 + b - 1}} \Gamma(b) \Gamma(b+\alpha_1 + \alpha_2-1)  \\
 \times \bigg ( {1 \over \Gamma(\alpha_1)\Gamma(b+\alpha_2) } {}_2 F_1 \Big (b , b+ \alpha_1 + \alpha_2-1;b+\alpha_2; {\beta_2 \over \beta_1 + \beta_2} \Big ) + (\alpha_1 \leftrightarrow \alpha_2, \: \beta_1 \leftrightarrow \beta_2 ) \bigg ).
\end{multline}
\end{prop}

Proofs of the above results are given in Section \ref{S2}, with a following discussion  in Section \ref{S3}.

\section{Proofs}\label{S2}
\subsection{Proof of Proposition \ref{P1} and Corollary \ref{C1}}
Taking the logarithmic derivative of (\ref{1.2}) gives
\begin{equation}\label{2.0a}
(1 - i t/\beta_1) (1 + i t/\beta_2)
\phi_{X_1 - X_2}'(t) = \Big ( (i \alpha_1/\beta_1)(1+it/\beta_2) -
(i \alpha_2/\beta_2)(1-it/\beta_1) \Big ) \phi_{X_1 - X_2}(t).
\end{equation}
Now multiply both sides by ${1 \over 2 \pi} e^{-ixt}$ and integrate over $t \in \mathbb R$.  We see that with the assumption $\alpha_1 + \alpha_2 > 2$ the integrand decays sufficiently fast that the integrals on both sides exist.

Focusing attention on the LHS, use of integration by parts gives
\begin{multline}\label{10.1}
{1 \over 2 \pi} \int_{-\infty}^\infty {\Big (   (i /\beta_1) -
(i / \beta_2) - (2t /\beta_1 \beta_2) \Big ) e^{-i x t} \over 
(1 - i t/\beta_1)^{\alpha_1} 
(1 + i t/\beta_1)^{\alpha_2}} \, dt +
{i x \over 2 \pi} 
\int_{-\infty}^\infty 
{(1 - i t/\beta_1) (1 + i t/\beta_2)e^{-i x t} \over
(1 - i t/\beta_1)^{\alpha_1} 
(1 + i t/\beta_1)^{\alpha_2}} \, dt \\
=  \Big ( {i \over \beta_1} -
{i \over \beta_2}\Big ) p_{X_1 - X_2}(x) + {2 \over i \beta_1 \beta_2}
p_{X_1 - X_2}'(x) + i x p_{X_1 - X_2}(x) + i x \Big ( {1 \over \beta_1} -
{1 \over \beta_2}\Big ) p_{X_1 - X_2}'(x) \\
 + {i x \over \beta_1 \beta_2}
\int_{-\infty}^\infty {t^2 e^{-i x t} \over (1 - i t/\beta_1)^{\alpha_1} 
(1 + i t/\beta_1)^{\alpha_2}} \, dt.
\end{multline}
 Here the equality follows from (\ref{1.3}) and the fact that for $\alpha_1 + \alpha_2 > 2$ the implied decay of the integrand allows for one differentiation under the integral sign, since the resulting integral is absolutely convergent. Strengthening the assumption to $\alpha_1 + \alpha_2 > 3$ allows the final term in (\ref{10.1}) to be identified as
 \begin{equation}\label{10.1a}
 - {i x \over \beta_1 \beta_2}
 p_X''(x).
 \end{equation}

 On the RHS, the assumption $\alpha_1 + \alpha_2 > 2$ is sufficient for the resulting integral to be identified as 
 \begin{equation}\label{10.1b}
 \Big ( {i \alpha_1 \over \beta_1} -
{i \alpha_2 \over \beta_2}\Big ) p_{X_1 - X_2}(x) + 
{1 \over i \beta_1 \beta_2} (\alpha_1 + \alpha_2) p_{X_1 - X_2}'(x).
 \end{equation}
 Substituting (\ref{10.1a}) for the final term in (\ref{10.1}) and equating with (\ref{10.1b}) gives the differential equation (\ref{2.7}), proved at this stage under the assumption $\alpha_1 + \alpha_2 > 3$.
 To extend the validity of (\ref{2.7}) to all $\alpha_1, \alpha_2 > 0$ we examine the dependency on the variables of the various terms as functions of complex $\alpha_1, \alpha_2$ in the domain ${\rm Re} \, \alpha_1 > 0$ and ${\rm Re} \, \alpha_2 > 0$.

 First, we see from the large $t$ behaviour of the integrand in (\ref{1.3}) that $p_{X_1 - X_2}(x)$ is analytic in $\alpha_1, \alpha_2$ 
 for ${\rm Re} \, \alpha_1 > 0$ and ${\rm Re} \, \alpha_2 > 0$
 for all fixed $x$ excluding $x=0$ (to include $x=0$ would require the extra condition ${\rm Re} \, (\alpha_1 + \alpha_2) > 1$). On this one sees the validity of repeated differentiation with respect to $\alpha_1, \alpha_2$ under the integral sign, justified by the absolute convergence of the resulting integral. However we cannot immediately make the same conclusion in relation to $p_{X_1 - X_2}'(x), p_{X_1 - X_2}''(x)$ as only with the extra conditions ${\rm Re} \, (\alpha_1 + \alpha_2) > 2$ and ${\rm Re} \, (\alpha_1 + \alpha_2) > 3$ respectively can justification of taking the differentiation inside of the integrand be made using the criteria that the integrals are then absolutely convergent.

 Fortunately, following an idea in \cite[Paragraph above Remark 2.10]{Ga21} (see also \cite[\S 2.1]{ZKKK20}) functional forms of $p_{X_1 - X_2}'(x), p_{X_1 - X_2}''(x)$ for general ${\rm Re} \, \alpha_1 > 0$ and ${\rm Re} \, \alpha_2 > 0$ can be given, and from these forms the sought analyticity can be concluded. Thus in the integrand of (\ref{1.3}) the exponential function $e^{-ixt}$ is written in terms of a derivative with respect to $t$, and then integration by parts is carried out. This gives
 the identity
 \begin{equation}\label{PpP}
 p_{X_1 - X_2}(x) = {1 \over x} \bigg (
 \Big ( {\alpha_1 \over \beta_1} - 
 {\alpha_2 \over \beta_2} \Big )
 p_{X_1 - X_2}(x) \Big |_{+,+} -
\Big ( {\alpha_1 + \alpha_2 \over \beta_1 \beta_2} \Big )
 p_{X_1 - X_2}'(x) \Big |_{+,+} \bigg ),
 \end{equation}
 where the notation $|_{+,+}$ indicates that the parameters $\alpha_1, \alpha_2$ are each to be increased by 1. In this expression the derivative term $p_{X_1 - X_2}'(x) |_{+,+}$ permits 
 (is defined by) 
 differentiation under the integral sign, which is justified by the absolute convergence of the integral. Furthermore, if we iterate this identity one more time, we see that the decay of all the resulting integrands on the RHS is, for  
${\rm Re} \, \alpha_1 > 0$ and ${\rm Re} \, \alpha_2 > 0$, faster than $1/|t|^2$ for large $t$. Hence  differentiation with respect to $x$ can be carried out under the integrand in this parameter range, giving rise to integrals which by inspection are analytic for ${\rm Re} \, \alpha_1 > 0$ and ${\rm Re} \, \alpha_2 > 0$. This then holds true of $p_{X_1 - X_2}'(x)$, i.e.~the derivative with respect to $x$ of the LHS of (\ref{PpP}), except possibly for $x=0$ as seen by the factor of $1/x$ on the RHS of (\ref{PpP}). A further iteration of (\ref{PpP}) (bringing the total to three) allows for the same conclusion in relation to $p_{X_1 - X_2}''(x)$.

Thus we have the situation that the identity (differential equation) (\ref{2.7}) between functions analytic in ${\rm Re} \, \alpha_1 > 0$ and ${\rm Re} \, \alpha_2 > 0$ has proved in the dense subset of this domain ${\rm Re} \, ( \alpha_1 + \alpha_2) > 3 $.  The uniqueness of analytic continuation extends the identity to all of the domain of analyticity.

With (\ref{2.7}) established, multiplying through by $g(x)$ and integrating by parts gives (\ref{1.18}). Here it assumed that the behaviour of $g(x)$ for large $|x|$ is such that the implied expectations are well defined. Choosing in (\ref{1.18}) $g(x) = x^{k}$ for $k$ a positive integer gives (\ref{1.19}).

\subsection{Proof of Propositions \ref{P3} and \ref{P4} }
In relation to Proposition \ref{P3}, we return to the definition of the gamma difference random variable $X$ as $X = X_1 - X_2$, where $X_1, X_2$ are gamma random variables from the distributions $\Gamma[\alpha_1,\beta_1]$ and $\Gamma[\alpha_2,\beta_2]$ respectively. Thus, as noted in \cite[second displayed equation in Section 3]{Kl15} it follows
\begin{equation}\label{2.5i}
\mathbb E( (X_1 - X_2)^k ) =
\sum_{l=0}^k \binom{k}{l} (-1)^l 
\mathbb E( X_1^{k-l}) \mathbb E( X_2^{l}).
\end{equation}

For $X_2$ a gamma random variable from $\Gamma[\alpha_2,\beta_2]$, we can compute directly from the corresponding probability density function (\ref{1.1}) that for $l$ a non-negative integer
\begin{equation}\label{2.6i}
\mathbb E(X_2^l) = \beta_2^{-l} {\Gamma(\alpha_2+l) \over \Gamma(\alpha_2)} =
\beta_2^{-l} (\alpha_2)_l,
\end{equation}
where $(\alpha_2)_l$ denotes the increasing Pochhammer symbol. Starting with (\ref{2.6i}) as it applies to $\mathbb E(X_1^l)$, then replacing $l$ by $k-l$ and 
manipulating $(\alpha_1)_{k-l}$
shows
\begin{equation}\label{2.7i}
\mathbb E(X_1^{k-l}) = (-1)^l \beta_1^{l-k} {\prod_{i=1}^k (\alpha_1 + i -1) \over (-(\alpha_1 + k - 1))_l}. 
\end{equation}
Furthermore, we note that in terms of the Pochhammer symbol, we can write
\begin{equation}\label{2.8i}
\binom{k}{l} = {(-1)^l (-k)_l \over l!}.
\end{equation}

Substituting (\ref{2.6i}), (\ref{2.7i}) and (\ref{2.8i}) in (\ref{2.5i}) and recalling the series form of the degree $k$ hypergeometric polynomial
${}_2 F_1(-k,b;c;z)$,
\begin{equation}\label{2.9i}
{}_2 F_1(-k,b;c;z) = \sum_{l=0}^k {(-k)_l (b)_l \over l! (c)_l} z^l,
\end{equation}
the first of the functional forms in (\ref{1.20}) results. To obtain the second expression, we note from the fact $(X_1 - X_2)^k = (-1)^k (X_2 - X_1)^k$ that the only effect of the interchange $(\alpha_1, \beta_1) \leftrightarrow (\alpha_2, \beta_2)$ is a factor of $(-1)^k$.

We turn our attention now to the derivation of (\ref{1.21}) in Proposition \ref{P4}. Our starting point is (\ref{1.6}), from which we read off
\begin{multline}\label{2.11x}
\mathbb E(|X|^{b - 1} ) =
{\beta_1^{\alpha_1} \beta_2^{\alpha_2} \over (\beta_1 + \beta_2)^{\alpha_1 + \alpha_2 + b - 1}} \\
\times \bigg ( {1 \over \Gamma(\alpha_1)}  \int_0^\infty 
  x^{b-1} e^{-\beta_1 x/(\beta_1 + \beta_2)} 
 U(1-\alpha_1,2-\alpha_1 - \alpha_2; x ) \, dx + (\alpha_1 \leftrightarrow \alpha_2, \: \beta_1 \leftrightarrow \beta_2 ) \bigg ).
 \end{multline}
 The integral in this expression is tabulated in \cite[Eq.~7.621(6)]{GR80}, where subject to the requirements that $b>{\rm max}\, (0, 1 - \alpha_1- \alpha_2)$ it is expressed in terms of a particular ${}_2 F_1$ hypergeometric function, implying (\ref{1.21}).

\section{Discussion}\label{S3}
It is a classical result that the Tricomi function $U(a,b;z)$ satisfies the second order differential equation known as Kummer's equation \cite[Eq.~13.2.1]{DLMF}. Combining this knowledge with the functional form (\ref{1.6}) allows for the differential equation characterisation of $p_{X_1 - X_2}(x)$ (\ref{2.7}) to be independently verified.

If we substitute $g(x) = e^{itx}$ in (\ref{1.18}) we recover the first order equation (\ref{2.0a}) for the characteristic function $\phi_{X_1 - X_2}(t)$. This shows that validity of (\ref{1.18}) for general twice differentiable $g(x)$ implies that the random variable $X$ is uniquely determined as having gamma difference distribution.

The structure of (\ref{2.0a}) can be generalised to the class of characteristic function $\phi_X(t)$ say satisfying the first order differential equation
$$
A_d(it) \phi_X'(t) = i B_{\tilde{d}}(it) \phi_X(t),
$$
for polynomials $A_d$ and $B_{\tilde{d}}$ of degrees $d$ and $\tilde{d}$ respectively. From this starting point, the Stein characterisation of the law of $X$,
$$
\mathbb E \bigg ( X A_d\Big ( {d \over dX} \Big ) g(X) - B_{\tilde{d}}\Big ( {d \over dX} \Big ) g(X)\bigg ) = 0,
$$
valid for a suitable class of $g(x)$, has been established in \cite[Lemma 2.1]{AAPS20}. Thus our (\ref{1.18}) is a special case of this more general result.\footnote{I thank R.E.~Gaunt for pointing this out to me.}

In relation to the recurrence (\ref{1.19}) for the moments, we see that it exhibits $m_k(X)$ to be a polynomial of degree $k$ in each of $1/\beta_1, 1/\beta_2,\alpha_1$ and $\alpha_2$. This structural feature is consistent with the exact evaluation (\ref{1.20}).

The first negative moment is of interest from the viewpoint an application to electro-optical imaging sensors \cite{He19}. While this quantity is ill-defined due to $p_{X_1 - X_2}(0) \ne 0$, for $\alpha_1 + \alpha_2 > 1$ it is
well defined in a principal value (PV) sense, where  according to \cite[proof of Th.~3.3]{He19} one has
$$
{\rm PV} \, \mathbb E(X^{-1}) = {\rm Re} \, \lim_{\epsilon \to 0} \Big ( e^{\pi i \epsilon} \mathbb E(|X|^{\epsilon - 1} \chi_{X \le 0}) +  \mathbb E(X^{\epsilon - 1} \chi_{X > 0}) \Big ).
$$
Each term on the RHS can be evaluated using the same integral evaluation  as used in (\ref{2.11x}), leading to an evaluation formula for ${\rm PV} \, \mathbb E(X^{-1})$ in terms of
particular ${}_3 F_2$ and digamma functions \cite[Th.~3.3]{He19}.

It must be that with $b-1 = k$, $k$ even, the positive integer moment formula (\ref{1.20}) agrees with the continuous moment formula (\ref{1.21}). To get some insight into this, we first note that the argument of the second (implied) ${}_2 F_1$ function in (\ref{1.21}) relates to the argument of the first ${}_2 F_1$ function by $z \mapsto 1 - z$. This suggest applying the transformation identity between ${}_2 F_1$ of argument $z$ and a linear combination of two ${}_2 F_1$ functions of argument $1-z$ \cite[Eq.~9.131(2)]{GR80}. Doing this gives the rewrite of (\ref{1.21})
\begin{multline} \label{1.21+} 
\mathbb E(|X|^{b-1}) = {\beta_1^{\alpha_1} \beta_2^{\alpha_2} \over (\beta_1 + \beta_2)^{\alpha_1 + \alpha_2 + b - 1}} \Gamma(b) \Gamma(b+\alpha_1 + \alpha_2-1)  \\
 \times \bigg ( {1 \over \Gamma(\alpha_1)\Gamma(b+\alpha_2) }
 \Big ( 1 + {\sin \pi \alpha_2 \over \sin \pi (b + \alpha_2)} \Big ) {}_2 F_1 \Big (b , b+ \alpha_1 + \alpha_2-1;b+\alpha_2; {\beta_2 \over \beta_1 + \beta_2} \Big )  \\ +
 {\Gamma(b + \alpha_2 - 1) \over \Gamma(\alpha_2) \Gamma(b) \Gamma(b+\alpha_1 + \alpha_2-1) } \Big ( {\beta_2 \over \beta_1 + \beta_2} \Big )^{-b-\alpha_2 + 1}
 {}_2 F_1 \Big (\alpha_1 ,  -\alpha_2 + 1; -b - \alpha_2+2; {\beta_2 \over \beta_1 + \beta_2} \Big )
 \bigg ).
\end{multline}
We see from this a simplification in the case that $b-1 = k$, $k$ even, namely that the prefactor of the first ${}_2 F_1$ function then vanishes, leaving us with
\begin{equation} \label{1.21+1}
\mathbb E(X^k) = {\beta_1^{\alpha_1} \beta_2^{-k} \over (\beta_1 + \beta_2)^{\alpha_1}} {\Gamma(k+\alpha_2) \over \Gamma(\alpha_2)} 
{}_2 F_1 \Big (\alpha_1 ,  -\alpha_2 + 1; -k - \alpha_2+1; {\beta_2 \over \beta_1 + \beta_2} \Big ).
\end{equation}
A standard Pfaff transformation mapping the ${}_2 F_1$ function with argument $z$ to a ${}_2 F_1$ function with argument $z/(1-z)$ shows that (\ref{1.21+1}) is equivalent to the second of the forms given in (\ref{1.20}).

We now specialise the parameters of the gamma difference distribution to those of the variance gamma distribution as specified above (\ref{1.2a}). Doing this in (\ref{2.7}) we see that the known differential equation (\ref{VG3}) for the probability density function of the variance gamma distribution is reclaimed. Similarly, we see that the three term recurrence for the positive integer moments (\ref{1.19}) reduces to the three term recurrence for the positive integer moments (\ref{2.6a+}) known for the variance gamma distribution. 
In relation to the formula 
(\ref{2.6b}) for the continuous moments, we apply to this the particular quadratic transformation formula \cite[Eq.~15.8.7]{DLMF} as applies to ${}_2 F_1 (a,b;1/2,z)$. We obtain the appropriate specialisation of the gamma difference result (\ref{1.21}).

Finally we comment on the appearance of the case $\beta_1 = \beta_2 = 1$, $\alpha_2 = \bar{\alpha}_1$ of (\ref{1.2}) --- to be denoted $w(t;\alpha_1)$ --- in random matrix theory, and an associated linear differential  and difference equation. In this field, one encounters the eigenvalue probability density function on the space of Hermitian matrices proportional to 
$$
\prod_{l=1}^N w(x_l;N+\alpha) \prod_{1 \le j < k \le N} | x_k - x_j|^2;
$$
see \cite[Eq.~(3.124)]{Fo10}. For this to be normalisable, it is required that ${\rm Re} \, \alpha > 1/2$.
This functional form is well known to give rise to a determinantal point process (see \cite[Ch.~5]{Fo10}), for which the eigenvalue density $\rho_{(1),N}(x)$  --- this quantity is normalised to integrate to $N$ --- has the explicit functional form
\begin{equation}\label{5.1}
\rho_{(1),N}(x) =  w(x;N+\alpha) \sum_{k=0}^{N-1} {1 \over h_k} (I_k^{(N+\alpha)}(x))^2.
\end{equation}
Here $\{ I_k^{(N+\alpha)}(x) \}_{k=0}^{N-1}$ are the so-called Romanovski polynomials, defined as the monic polynomials of degree $k$ exhibiting the orthogonality
\begin{equation}\label{5.2}
\int_{-\infty}^\infty w(x;N+\alpha) 
I_j^{(N+\alpha)}(x) I_k^{(N+\alpha)}(x) \, dx =
h_j \delta_{j,k}
\end{equation}
for $0 \le j,k \le N-1$.
See \cite[Eq.~(2.2)]{FLT21} in relation to (\ref{5.1}), \cite[Eq.~(2.4)]{FLT21} for an explicit hypergeometric polynomial expression for the $I_k^{(N+\alpha)}(x)$, and \cite[Eq.~(2.5)]{FLT21} for a gamma function evaluation of the normalisation $h_j$.

There are two points of interest in relation to $\rho_{(1),N}(x)$ in the present context. One is its characterisation as the solution of a third order linear differential equation valid for general positive integer $N$ \cite[Prop.~3.5]{FR21} (As we know from (\ref{2.0a}), the special case $N=1$ admits a simpler first order differential equation characterisation.) The other is that in the so-called symmetric case of $\alpha$ real, the sum of successive absolute continuous moments $\mu_k := m_{2k+2} + m_{2k}$, with $m_{2k} := \int_{-\infty}^\infty |x|^{2k} \rho_{(1),N}(x) \, dx$ (these are well defined for $-1/2 < {\rm Re} \, k < \alpha - 1/2$) satisfy a  three term recurrence \cite[Cor.~4.3]{FR21}. Moreover, the $\mu_k$ admit an explicit evaluation in terms of particular continuous Hahn polynomials \cite{ABGS21}.

\subsection*{Acknowledgements}
This work is part of a research program supported by the Australian Research Council  grant DP210102887.
Helpful remarks by R.E.~Gaunt on the first draft of this work are appreciated. I thank too A.~Hendrickson for drawing
my attention to \cite{He19}.

\providecommand{\bysame}{\leavevmode\hbox to3em{\hrulefill}\thinspace}
\providecommand{\MR}{\relax\ifhmode\unskip\space\fi MR }
\providecommand{\MRhref}[2]{%
  \href{http://www.ams.org/mathscinet-getitem?mr=#1}{#2}
}
\providecommand{\href}[2]{#2}

\end{document}